\documentclass[11pt]{article}
\newcommand{\be}{\begin{equation}}
\newcommand{\ee}{\end{equation}}
\newcommand{\ba}{\begin{array}}
\newcommand{\ea}{\end{array}}
\newcommand{\bea}{\begin{eqnarray}}
\newcommand{\eea}{\end{eqnarray}}
\newcommand{\bee}{\begin{eqnarray*}}
\newcommand{\eee}{\end{eqnarray*}}

\newtheorem{theorem}{Theorem}
\newtheorem{lemma}{Lemma}
\newtheorem{proposition}{Proposition}

\newtheorem{corollary}{Corollary}

\def\R{{\bf R}}

\def\lim{\mathop{\rm lim}}
\def\goto{\rightarrow}

\def\sup{\mathop{\rm sup}}
\def\exp{{\rm exp}}
\def\e{\varepsilon}

\def\l{\lambda}
\def\ch{{\rm ch}}

\def\ds{\displaystyle}

\def\ut{\widetilde u}
\def\et{\widetilde \e}
\def\ct{\widetilde c}
\def\xt{\widetilde x}

\def\norm#1{\left\|#1\right\|}
\def\eqref#1{(\ref{#1})}
\def\bket#1{\left\{#1\right\}}
\def\bkt#1{\left[#1\right]}
\def\bke#1{\left(#1\right)}

\def\k{\kappa}

\begin{document}
\title{Stability and asymptotic stability in the energy space \\
of the sum of
$N$ solitons  for subcritical gKdV equations }

\author{Yvan Martel, Frank Merle and Tai--Peng Tsai}
\date{}

\maketitle
\begin{abstract}
We prove in this paper the stability and asymptotic stability in 
$H^1$ of a decoupled sum of $N$ solitons for the subcritical generalized 
KdV equations $u_t+(u_{xx}+u^p)_x=0$ ($1<p<5$).
The proof of the stability result is based on energy arguments and monotonicity
of local $L^2$ norm. Note that the result is new even for $p=2$ 
(the KdV equation). The asymptotic stability result then follows directly
from a rigidity theorem in \cite{MM3}.
\end{abstract}

\section{Introduction}

In this paper, we consider the
generalized Korteweg--de Vries equations
\begin{equation}
\left\{
\begin{array}{l}
u_t+(u_{xx}+u^p)_x=0 , \quad (t,x)\in \R \times \R,\\[5pt]
u(0,x)= u_0(x) ,\quad x\in \R,
\end{array}
\right.\label{kdvp}
\end{equation}
for $1<p<5$ and $u_0\in H^1(\R)$.
This model for $p=2$
was first introduced in the study of waves on shallow water,
see Korteweg and de Vries \cite{KDV}. It also appears for $p=2$ and 3,
in other areas of Physics (see e.g. Lamb \cite{Lamb}).

\medskip

Recall that (\ref{kdvp}) is well--posed in the energy space $H^1$. 
For  $p=2,3,4$,
it was proved by Kenig, Ponce and Vega \cite{KPV} 
(see also Kato \cite{K}, Ginibre and Tsutsumi \cite{GT}),
that for
$u_0\in H^1(\R)$, there exists a unique  solution
$u\in C(\R,H^1(\R))$ of (\ref{kdvp}) satisfying the following two conservation 
laws, for all $t\in \R$,
\begin{equation}\label{masse}
  \int u^2(t)=\int u^2_0,
\end{equation}
\begin{equation}\label{energie}
  E(u(t))={1\over 2} \int u^2_x(t) -{1\over p+1} \int u^{p+1}(t)
         ={1\over 2} \int u_{0x}^2 -{1\over p+1} \int u_0^{p+1}.
\end{equation}
For $p=2,3,4$, global existence of all solutions
in $H^1$, as well as uniform bound in $H^1$,
follow directly from the Gagliardo--Nirenberg inequality,
$$\forall v\in H^1(\R),\quad 
\int |v|^{p+1} \le C(p)
\left( \int v^2 \right)^{p+3\over 4} \left(\int v_x^2\right)^{{p-1}\over 4},
$$
and relations (\ref{masse}), (\ref{energie}), giving
a uniform bound in $H^1$ for any solution.

This is in contrast with the case $p=5$, for which 
there exist solutions $u(t)$ of (\ref{kdvp})
such that $|u(t)|_{H^1}\goto +\infty$ as $t\goto T$, for
$0<T<+\infty$, see \cite{M3} and \cite{MM5}.
For $p>5$ such behavior is also conjectured. 
Thus, for the question of global existence and  bound in $H^1$,
the case $1<p<5$ is called the 
subcritical case, $p=5$ the critical case and $p>5$ the supercritical case.

\medskip

Equation (\ref{kdvp}) has explicit traveling wave solutions, called solitons,
which play a fundamental role in the generic behavior of the solutions.
Let
\be\label{ground}
Q(x)=\left( {p+1 \over 2\, \ch^2\left(
{p-1\over 2} x\right)}\right)^{1\over p-1}
\ee
be the only positive solution in $H^1(\R)$ (up to translation) of
$Q_{xx}+Q^p=Q,$
 and for $c>0$, let
$Q_c(x)=c^{1\over p-1} Q\left(\sqrt{c} x \right)$.
The traveling waves solutions of (\ref{kdvp}) are
$$u(t,x)
=Q_c(x- c t)=c^{1\over p-1} Q\left(\sqrt{c} ( x -ct)\right),$$
where $c>0$ is the speed of the soliton.

For the KdV equation ($p=2$),
there is a much wider class of special explicit solutions for (\ref{kdvp}),
called $N$--solitons. They correspond to the superposition of $N$
traveling waves  with different speeds
that interact and then remain unchanged
after interaction.
The $N$--solitons  behave asymptotically in large
time as the sum of $N$ traveling waves, and as for the single solitons,
there is no dispersion.
We refer to \cite{MIURA} for  explicit expressions and
further  properties  of these solutions.
For $p \not =2$,  even the existence of solutions behaving
asymptotically as the sum of $N$ solitons was not known.

\medskip

Important notions for these solutions are the stability and asymptotic stability
with respect to initial data.

{\sl For $c>0$, the soliton
$Q_c(x-ct)$ is stable in $H^1$ if:}
$$\hbox{
$\forall \delta_0>0, $ $ \exists \alpha_0>0$ $/$
$|u_0-Q_{c}|_{H^1}\le \alpha_0$  $\Rightarrow$
$
\forall t\ge 0,$ $\exists x(t)$ $/$ $|u(t)-Q_{c}(.-x(t))|_{H^1}\le \delta_0.
$
}$$

{\sl The family of solitons
$\{Q_c(x-x_0-c t),~c>0,~x_0\in \R \}$ is asymptotically stable if:}
$$\hbox{$
\exists \alpha_0>0$ $/$ $|u_0-Q_{c}|_{H^1}\le \alpha_0$ $\Rightarrow$
$
\exists c_{+\infty}, x(t)$ $/$ 
$\ds u(t,.+x(t))\mathop{\rightharpoonup}_{t\goto +\infty} Q_{c_{+\infty}}$
in $ H^1 $.}
$$

We recall  previously known results concerning
the notions of stability of solitons and $N$ solitons:

\smallskip

- In the subcritical case: $p=2,3,4$, it follows
from energetic arguments that the solitons are $H^1$ stable
(see Benjamin \cite{Be} and Weinstein \cite{W2}). Moreover, 
Martel and Merle \cite{MM3} prove the asymptotic stability
of the family of solitons in the energy space. The proof relies on
a rigidity theorem close to the family
of solitons, which was first given for the
critical case (\cite{MM2}), and which is based on nonlinear argument.
(Pego and Weinstein \cite{PW} prove this result
for $p=2,3$ for initial data with exponential decay as $x\goto +\infty$.)

In the case of the KdV equation,
Maddocks and Sachs \cite{MS}
prove the stability in $H^N(\R)$ of  $N$--solitons
(recall that there are explicit solutions of the KdV equation) :
for any  initial data $u_0$ close in $H^N(\R)$ to an $N$--soliton,
the solution $u(t)$ of the KdV equation remains uniformly close in $H^N(\R)$
for all time to an $N$ soliton profile with same speeds.
Their proof involves $N$ conserved quantities for the KdV equation,
and this is the reason why they need to impose closeness in
high regularity spaces. Note that this result is known 
only with $p=2$ and with this regularity assumption of the initial data.
Asymptotic stability is unknown in this context.

\smallskip
- In the critical case $p=5$,
any solution with negative energy
initially close to the soliton  blows up in finite
or infinite time in $H^1$ (Merle \cite{M3}), and actually 
blows up in finite time
if the initial data satisfies in addition a
polynomial decay condition on the right in space (Martel and Merle \cite{MM5}).
(Note that $E(Q)=0$ for $p=5$.) Of course this implies the instability
of the soliton.
These results rely on rigidity theorems around the soliton.

\smallskip

- In the supercritical case $p>5$,
Bona, Souganidis, and
Strauss \cite{BSS} proved, using Grillakis, Shatah, and Strauss \cite{GSS} type
arguments, $H^1$ instability of solitons. Moreover,
numerical experiments, see e.g. Dix and McKinney
\cite{DK}, suggest existence of blow up solutions arbitrarily
close to the family of solitons.

\bigskip

In this paper, for $p=2,3,4$, using techniques developed for the 
critical and subcritical cases in \cite{MM2} and \cite{MM3}
as well as a direct variational argument in $H^1$,
we prove the stability and asymptotic stability
of the sum
$$\sum_{j=1}^N Q_{c_j^0}(x-x_j), \quad \hbox{where \quad
$0<c_1^0<\ldots<c_N^0,$ \quad  $x_1<\ldots<x_N$,}
$$
in $H^1(\R)$, for $t\ge 0$.

\begin{theorem}[Asymptotic stability of the sum of
$N$ solitons] \label{MAINTH}
Let $p=2$, 3 or 4. Let $0<c_1^0<\ldots<c_N^0.$
There exist $\gamma_0,A_0,L_0,\alpha_0>0$ such that the following
is true.
Let $u_0\in H^1(\R)$ and assume that there exist  $L>L_0$, $\alpha<\alpha_0$,
and $x_1^0<\ldots<x_N^0$, 
such that
\be\label{etoileee}
\Bigl| u_0-\sum_{j=1}^N Q_{c_j^0}(.-x_j^0) \Bigr|_{H^1}\le \alpha,
\quad \hbox{and $x_j^0>x_{j-1}^0+L$, for all $j=2,\ldots,N$}.
\ee
 Let $u(t)$ be the solution of (\ref{kdvp}).
Then, there exist $x_1(t),\ldots,x_N(t)$  such that

(i) Stability of the sum of $N$ decoupled solitons.
\quad 
\be\label{stabilityth}
\forall t\ge 0,\quad \Bigl| u(t) - \sum_{j=1}^N Q_{c_j^0} (x-x_j(t)) \Bigr|_{H^1}
\le A_0 \left(\alpha + e^{-\gamma_0 L}\right).
\ee

(ii) Asymptotic stability of the sum of $N$ solitons. 
Moreover, there exist $c_1^{+\infty},\ldots,
c_N^{+\infty}$, 
with 
$|c_j^{+\infty}-c_j^0|\le A_0 \left(\alpha + e^{-\gamma_0 L}
\right),
$
such that
\be\label{asympstabth}
 \Bigl| u(t) - \sum_{j=1}^N Q_{c_j^{
+\infty}} (x-x_j(t)) \Bigr|_{L^2(x>{c_1^0 t /10})}
\goto 0, \quad \dot x_j(t)\goto c_j^{+\infty}
 \quad  \hbox{as $t\goto +\infty$}.
\ee
\end{theorem}

{\bf Remark 1.}  \quad 
It is well-known that for $p=2$ and $p=3$, (\ref{kdvp})
is completely integrable. Indeed,
for suitable $u_0$ ($u_0$ and its derivatives with exponential
decay at infinity) there exist an infinite number of conservation laws,
see e.g. Lax \cite{LAX} and Miura \cite{MIURA}. Moreover, 
many results on these equations
rely on the inverse scattering method, which transform the problem in
a sequence of linear problems (but requires strong decay assumption on
the solution). 
{\sl In this paper, we do not use integrability.}

\medskip

{\bf Remark 2.} \quad  For  Schr\"odinger type equations,
Perelman \cite{P} and Buslaev and Perelman \cite{BP}, with strong 
conditions on initial data and nonlinearity, and using a linearization method
around the soliton, prove asymptotic stability results by a fixed point 
argument. Unfortunately, this method breaks down without decay assumption 
on the initial data. 

\medskip

{\bf Remark 3.} \quad In Theorem \ref{MAINTH} (ii), we cannot have
convergence to zero in $L^2(x>0)$. Indeed, assumption (\ref{etoileee})
on the initial data allows the existence in $u(t)$ of an additional
soliton of size less that $\alpha$ (thus traveling at arbitrarily
small speed). For $p=2$, an explicit example can be constructed 
using the $N$--soliton solutions.

\medskip

Recall that for $p=2$  any $N$-soliton solution
has the form
$v(t,x)=U^{(N)}(x;c_j,x_j-c_j t)$, where 
$\{U^{(N)}(x;c_j,y_j);c_j>0,y_j\in \R\}$ is the family of
explicit $N$-soliton profiles (see e.g. \cite{MS}, \S 3.1).
As a direct corollary of Theorem 1, for $p=2$, we prove 
stability and asymptotic stability of this family.

\begin{corollary}[Asymptotic stability in $H^1$
of $N$-solitons for $p=2$]
Let $p=2$. Let  $0<c_1^0<\ldots<c_N^0$ and $x_1^0,\ldots,x_N^0\in \R$.
For all $\delta_1>0$, there exists $\alpha_1>0$ such that the following
is true. Let $u(t)$ be a solution of (\ref{kdvp}).
If
$| u(0)- U^{(N)}(\,.\,;c^0_j,-x_j^0)|_{H^1}\le \alpha_1,$
 then there exist $x_j(t)$ such that
\be\label{stabilitycor}
\forall t>0,\quad | u(t) - U^{(N)}(\,.\,;c^0_j,-x_j(t))|_{H^1}
\le \delta_1.
\ee
Moreover,
there exist $c_j^{+\infty}>0$  such that
\be\label{asympstabcor}
 \Bigl| u(t) - U^{(N)}(\,.\,
;c_j^{+\infty},-x_j(t)) \Bigr|_{L^2(x>{c_1^0 t /10})}
\goto 0, \quad \dot x_j(t)\goto c_j^{+\infty}\quad  \hbox{as $t\goto +\infty$}.
\ee
\end{corollary}

Note that this  improves the result in  \cite{MS}
in two ways. First, stability is proved in $H^1$ instead of
$H^N$. Second, we also prove asymptotic stability as $t\goto +\infty$.
Corollary 1 is proved at the end of \S 4.

\medskip

Let us sketch the proof of these results.
For Theorem \ref{MAINTH}, using modulation theory, 
$u(t)=\sum_{j=1}^N Q_{c_j(t)}(x-x_j(t)) +\e(t,x)$,
where $\e(t)$ is small in $H^1$, and  $x_i(t)$, $c_i(t)$
are  geometrical parameters (see \S 2).
The stability result is equivalent to control
both the variation of $c_j(t)$ and the size of $\e(t)$ in $H^1$
(\S 3).

Our main arguments are based on $L^2$ properties of the solution.
From \cite{MM2} and \cite{MM3}, the $L^2$ norm 
of the solution at the right of each soliton is almost decreasing in time.
This property together with energy argument allows us to prove that
the  variation of $c_j(t)$ is quadratic
in $|\e(t)|_{H^1}$, which is a key of the problem.

Let us explain the argument formally by taking $\e=0$ and
so $u(t)=\sum Q_{c_j(t)} (x-x_j(t)).$ The energy
conservation becomes
$$\sum c_j^{\beta +{1\over 2}}(t) = \sum c_j^{\beta+{1\over 2}}(0),
$$
where $\beta={2\over p-1}$. The monotonicity of the $L^2$ norm
at the right of each soliton gives us
$$\Delta_j(t)=\sum_{k=j}^N c_k^{\beta-{1\over 2}}(t)
-c_k^{\beta-{1\over 2}}(0)\le 0.$$
We claim that $c_j(t)=c_j(0)$ by a convexity argument.
Indeed,
\bee 0 &=&\sum c_j^{\beta +{1\over 2}}(t) -  c_j^{\beta+{1\over 2}}(0)
\sim {2\beta+1\over 2\beta-1}
\sum c_j(0) (c_j^{\beta -{1\over 2}}(t) -  c_j^{\beta-{1\over 2}}(0))
\\
&=&{2\beta+1\over 2\beta-1} \sum
(c_j(0)-c_{j+1}(0))\Delta_j(t) \ge \sigma_0 \sum |\Delta_j(t)|
\ge \sigma_1 \sum|c_j(t)-c_j(0)|.\eee
Thus $c_j(t)$ is a constant at the first order. In fact,
we prove that  the variation in time
of $c_j(t)$ is of order 2 in $\e(t)$.

Then, we control the variation of $\e(t)$ in $H^1$ by a
refined version of this argument, using suitable 
orthogonality conditions on $\e$.

The asymptotic stability result
follows directly from a rigidity property of the flow of equation
(\ref{kdvp}) around the solitons (see \cite{MM3}) and 
monotonicity properties of the mass (\S 4).

\medskip

\noindent 
{\bf Acknowledgments.} \quad Part of this work was done
when Tai-Peng Tsai was visiting the University of Cergy--Pontoise,
whose hospitality is gratefully acknowledged.

\section{Decomposition and properties of a solution close to the sum
of $N$ solitons}

\subsection{Decomposition of the solution and conservation laws}
Fix $0<c_1^0<\ldots<c_N^0$  and
let 
$$
\sigma_0={1\over 2}\min(c_1^0,c_2^0-c_1^0,c_3^0-c_2^0,\ldots,c_N^0-c_{N-1}^0)
.$$
From modulation theory, we claim.

\begin{lemma}[Decomposition of the solution]\label{DECOMPOSITION}
 There exists $L_1, \alpha_1, K_1>0$ such that the
following is true. If for $L>L_1$, $0<\alpha<\alpha_1$,
$t_0>0$, we have
\be\label{close}
\sup_{0\le t\le t_0} \Big( \inf_{y_j>y_{j-1}+L}\Big\{
\Big|u(t,.)-\sum_{j=1}^N Q_{c^0_j}(.-y_j)\Big|_{H^1}
\Big\}\Big)<\alpha,
\ee
then there exist unique $C^1$ functions
$c_j:[0,t_0]\goto (0,+\infty),$ $x_j:[0,t_0]\goto \R$, such that
\be\label{defofeps}
\e(t,x)=u(t,x)-\sum_{j=1}^N R_j(t,x),
\quad \hbox{where} \quad R_{j}(t,x)=Q_{c_j(t)}(x-x_j(t)),
\ee
satisfies the following orthogonality conditions
\be\label{ortho}
\forall j, \forall t\in [0,t_0],\quad \int R_j(t) \e(t)=
\int (R_j(t))_x \e(t)=0.
\ee
Moreover,  there exists $K_1>0$ such that 
$\forall t\in [0,t_0],$
\be\label{smallness1}
|\e(t)|_{H^1}+
\sum_{j=1}^N |c_j(t)-c^0_j|
\le K_1 \alpha,
\ee
\be\label{faux13}
\forall j,~
\left|\dot c_j(t)\right|
+\left|\dot x_j(t)-c_j(t)\right|
\le K_1 \Big(\int e^{-\sqrt{\sigma_0} |x-x_j(t)|/2} \e^2(t)\Big)^{1/2} 
+ K_1 e^{-\sqrt{\sigma_0} (L+\sigma_0 t)/4}.
\ee
\end{lemma}

{\bf Proof.}\quad Lemma \ref{DECOMPOSITION} is a consequence of
Lemma \ref{MODULATION} (see Appendix) and standard arguments.
We refer to \cite{MM1} \S 2.3 for a complete proof in the case of a
single soliton.
In particular, 
$\e(t)$ satisfies $\forall t\in [0,t_0]$,
$$
\e_t + \e_{xxx} =
- \sum_{j=1}^N {\dot c_j\over 2 c_j} \left({2R_j\over p-1}+(x-x_i)(R_j)_x
\right)
 + \sum_{j=1}^N(\dot x_j-c_j) R_{jx}
-\Big(\Big(\e+\sum_{j=1}^N R_j\Big)^p - \sum_{j=1}^N R_j^p\Big)_x.
$$
By taking (formally) the scalar product of this equation by 
$R_j$ and $(R_j)_x$, and using calculations
in the proof of Lemma \ref{MODULATION}, we prove
$$|\dot c_j(t)|+|\dot x_j(t)-c_j(t)|
\le C\Big( \int e^{-\sqrt{\sigma_0} |x-x_j(t)|/2} \e^2(t)\Big)^{1/2}
+C \sum_{k\not = j} e^{-\sqrt{\sigma_0} |x_k(t)-x_j(t)|/2}.
$$
For $\alpha>0$ small enough, and $L$ large enough, we have
$|x_k(t)-x_j(t)|\ge {L\over 2}
+\sigma_0 t,$ and this proves (\ref{faux13}).

\bigskip

Next, by  using the conservation of energy for $u(t)$, i.e.
\[
E(u(t)):=\int \frac 12 \, u_x^2(t,x) - \frac 1{p+1} \, u^{p+1}(t,x) \, d x
=E(u_0),\]
and linearizing the energy around $R=\sum_{j=1}^N R_j$, we 
prove the following result.

\begin{lemma}[Energy bounds] \label{th:3}
There exist $K_2>0$ and $L_2>0$ such that the following
is true.
Assume that  $\forall j$, $ c_j(t)\ge \sigma_0,$ and 
$x_j(t)-x_{j-1}(t)\ge L\ge L_2.$
Then, $\forall t\in [0,t_0]$,
\bea && 
\left| 
\sum_{j=1}^N \bkt{E(R_j(t)) - E(R_j(0)) }
 +\frac 12 \int \, (\e_x^2 - p R^{p-1}\e^2)(t)  \right |
\nonumber \\ && \quad \le
K_2 \bket{ |\e(0)|_{H^1}^2 
+ |\e(t)|_{H^1}^3
+ e^{-\sqrt{\sigma_0} L/2}},   \label{29} 
\eea
where $K_2$ is a constant.
\end{lemma}

\noindent{\bf Proof.} \quad
Insert \eqref{defofeps} into $E(u(t))$ and integrate by parts. We have
\bea
E(u(t)) &=&\int \frac 12 \, R_x^2 - \frac 1{p+1}\, R^{p+1}
\, d x 
 -\int  \bke{R_{xx} + R^p} \e  \, d x \label{3}
+ \int  \frac 12 \, \e_x^2 - \frac p2 R^{p-1}\e^2 \, d x
\qquad  \label{4}
\\
&& + \int  \frac 1{p+1}\, \bke{-(R+\e)^{p+1} + R^{p+1}}
+ R^p \e+ \frac p2 R^{p-1}\e^2  \, d x   \label{5}
\eea
We first observe that
$
|\eqref{5}| \le C\norm{\e}_{H^1}^3.
$
Next, remark that 
$
\sigma_0\le c_j(t),$ $x_{j}(t) - x_{j-1}(t) \ge L,
$
implies 
$
|R_j(x,t)|+|(R_j)_x(x,t)|
\le C e^{- \sqrt{\sigma_0}|x-x_j(t)|}$, and so
\be \label{6}
\left|\int R_j(t) \, R_k(t) \, d x \right| + \left|\int (R_j)_x(t) \,
(R_k)_x(t) \, d x \right|\le C e^{- \sqrt{\sigma_0}L/2} 
\quad {\rm if} \ j \not = k.
\ee
Thus, by $(R_j)_{xx}+R_j^p=c_j R_j$, we have 
\be
\left| \eqref{4} - \sum_{j=1}^N E(R_j(t))+\int  \sum_j c_j R_j \e (t)
- \frac 12 \int \, (\e_x^2 - p R^{p-1}\e^2)(t)
\right|\le C e^{- \sqrt{\sigma_0}L/2}. 
\ee
From
$\int R_j(t) \e(t)=0$, we obtain
\bee 
 \left| E(u(t)) -\sum_{j=1}^N E(R_j(t)) 
- \frac 12 \int \, (\e_x^2 - p R^{p-1}\e^2)(t)\right |
 \le Ce^{- \sqrt{\sigma_0}L/2}+ C \norm{\e(t)}_{H^1}^3.
\eee
Since $E(u(t))=E(u(0))$, applying the previous formula at
$t=0$ and at $t$, we prove the lemma.


\subsection{Almost monotonicity of the mass at the right}
We follow the proof of Lemma 20 in \cite{MM2}.
Let 
\be \phi(x)=c Q(\sqrt{\sigma_0} x/2), \quad \psi(x)=\int_{-\infty}^x
\phi(y)dy, \quad {\rm where }\ c= \bke{{2\over \sqrt{ \sigma_0 } }
 \int_{-\infty}^\infty
Q}^{-1}. \ee Note that 
$\forall x\in
\R,$ $\psi'>0$,  $0< \psi(x)<1,$ and $\ds \lim_ {x\to -\infty}\psi(x)=0,$
$\ds \lim_ {x\to +\infty}\psi(x)=1.$
Let
\be j\ge 2, \quad
{\cal I}_j = \int u^2(t,x) \psi(x-m_j(t) )\, dx, \qquad
m_j(t)=\frac {x_{j-1}(t)+x_{j}(t)} 2.
\ee

\begin{lemma}[Almost monotonicity of the  mass on the right 
of each soliton \cite{MM2}] \label{th:2}
There exist $K_3=K_3(\sigma_0)>0$, $L_3=L_3(\sigma_0)>0$ such that 
the following is true. Let
 $t_1\in [0,t_0]$. Assume that
$\forall t\in [0,t_1],$ $\forall j$,
\be\label{al1}
\dot x_1(t)\ge \sigma_0, 
\quad 
\dot x_j(t)- \dot x_{j-1}(t)\ge \sigma_0,
\quad 
c_j(t)>\sigma_0,\quad \hbox{and} \quad 
|\e(t)|^{p-1}_{H^1}\le {\sigma_0\over 8\cdot 2^{p-1}}.
\ee
If for $L>L_3$, $\forall j\in \{2,\ldots,N\}$,
 $x_j(0)- x_{j-1}(0)\ge L$, then
\[
  {\cal I}_j(t_1) - {\cal I}_j(0)\le K_3\, e^{- \sqrt{\sigma_0} L /8}.
\]
\end{lemma}

\noindent{\bf Proof} \quad 
Let $j\in \{1,\ldots,N\}$.
Using equation \eqref{kdvp} and
integrating by parts several times, we have (see \cite{MM3}
equation (20)),
\[
\frac d{dt} {\cal I}_j(t)= \int \left( - 3 u_x^2 - \dot m u^2 +
\frac {2p}{p+1} u^{p+1} \right)\psi' + u^2 \psi^{(3)}.
\]

By definition of $\psi$,
$\psi^{(3)} \le \frac {\sigma_0}4 \psi'$, so that
\be\label{faux24}
\int u^2 \psi^{(3)}\le {\sigma_0 \over 4}\int u^2 \psi' .
\ee

To bound $\int u^{p+1}\psi'$, we divide the real line
to two regions: $I=[a,b]$ and its complement $I^C$, where $a=a(t)=
x_{j-1}(t)+\frac {L}4$ and $b=b(t)= x_{j}(t)-\frac {L}4$. 
Inside the interval $I$ we have
\[
\left| \int_I u^{p+1}\psi' \right|
\le \int u^2 \psi' \cdot \sup_I |u|^{p-1}
\]
Since for $x\in I$, for all $k=1,2,\ldots, N$,
$|x-x_k(t)|\ge {L\over 4},$ we have
\[
|u(t,x)|^{p-1}=\left|\sum _{k=1}^N R_k(t,x) +\e(t,x)\right|^{p-1}
\le C e^{-\sqrt{\sigma_0} L /4}  + 2^{p-1} |\e(t)|_{L^\infty}^{p-1}
\le  \frac {\sigma_0} 4,
\]
for $L>L_3(\sigma_0)$. Thus,
\be\label{vend1}
\left| \int_I u^{p+1}\psi' \right| \le  
 \frac {\sigma_0} 4 \int u^2 \psi'.
\ee
Next, in $I^C$, by Gagliardo Nirenberg inequality,
\bea
\int_{I^C} u^{p+1}\psi' \, dx & \le & \int u^{p+1} \, dx \cdot
\sup_{I^C} \psi' 
\le C \norm{u}_{H^1}^{p+1} \cdot  \exp \bket{-\frac {\sqrt{\sigma_0}}4
 \,[x_{j}(t)-x_{j-1}(t)-{\textstyle {L\over 2}}]}
\nonumber \\ &\le& C e^{-\frac {\sqrt{\sigma_0}}8
 (2\sigma_0 t +L)},\label{vend2}
\eea
by $x_j(t)-x_{j-1}(t)\ge x_j(0)-x_{j-1}(0)+\sigma_0 t
\ge L+\sigma_0 t$.
From $\dot m\ge \sigma_0$, \eqref{faux24},
\eqref{vend1} and \eqref{vend2}, we obtain
\[
\frac d{dt} {\cal I}_j(t) \le \int \left( - 3 u_x^2 - {\sigma_0
\over 2} u^2  \right)\psi'\, dx +C e^{-{\sqrt{\sigma_0}\over 8}
 (2\sigma_0 t +L)}
\le C e^{-{\sqrt{\sigma_0}\over 8} (2\sigma_0 t +L)} .
\]
Thus, by integrating between $0$ and $t_1$, we obtain
the conclusion. Note that $K_3$ and $L_3$ are chosen independently
of $t_1$.

\subsection{Positivity of the quadratic form}

By the choice of orthogonality conditions on $\e(t)$
and standard arguments, we claim the following lemma.

\begin{lemma}[Positivity of the quadratic form]\label{QUAD}
There exists $L_4>0$ and $\l_0>0$ such that if $\forall j$,
$c_j(t)\ge \sigma_0$,
$x_j(t)\ge x_{j-1}(t)+L_4$ then,
$\forall t\in [0,t_0]$,
\be\label{aprouver}
\int \e^2_x (t) - p R^{p-1}(t) \e^2(t) + c(t,x) \e^2(t) \ge 
\l_0 |\e(t)|^2_{H^1},
\ee
where $c(t,x)=c_1(t)+\sum_{j=2}^N (c_j(t)-c_{j-1}(t)) \psi(x-m_j(t)).
$
\end{lemma}

{\bf Proof of Lemma \ref{QUAD}.}
\quad 
It is well known 
that there exists $\l_1>0$ such that
if $v\in H^1(\R)$ satisfies $\int  Q v=\int Q_x v=0$, then
\be\label{single}
\int v_x^2-p Q^{p-1} v^2 + v^2\ge \l_1 |v|_{H^1}^2.
\ee
(See proof of Proposition 2.9 in Weinstein \cite{W3}.)
Now we give a local version of (\ref{single}).
Let $\Phi\in {\cal C}^2(\R)$, $\Phi(x)=\Phi(-x)$, $\Phi'\le 0$ on $\R^+$, 
with
$$\Phi(x)=1 \hbox{ on $[0,1]$;}\quad 
\Phi(x)=e^{-x} \hbox{ on $[2,+\infty)$,}\quad
e^{-x}\le \Phi(x)\le 3 e^{-x} \quad \hbox{on
$\R^+$.}$$ 
Let $\Phi_{B}(x)=\Phi\left({x\over B}\right).$
The following claim is similar to a part of the proof of some local Virial
relation in \S 2.2 of \cite {MM4}; see 
Appendix A, Steps 1 and 2, in \cite{MM4} for its proof.

\medskip

{\sl Claim.}\quad There exists $B_0>0$ such that, for all $B>B_0$,
if $v\in H^1(\R)$ satisfies
$\int  Q v=\int Q_x v=0$, then
\be\label{singlelocal}
\int \Phi_{B} \left(v_x^2-p Q^{p-1} v^2 + v^2\right)
\ge {\l_1\over 4} \int \Phi_{B} (v_x^2+v^2).
\ee

\medskip

We finish the proof of Lemma \ref{QUAD}.
Let $B>B_0$ 
to be chosen later 
and $L_4=4 k B$, where $k>1$ integer is to be chosen later.
We have 
\bee
\int \e^2_x  - p R^{p-1} \e^2 + c(t,x) \e^2 
& = & 
 \sum_{j=1}^N \int \Phi_{B}(x-x_j(t))
\left(\e^2_x - p R_j^{p-1} \e^2 + c_j(t)\e^2\right) 
\\ && 
-p \int \Big(R^{p-1}-
\sum_{j=1}^N  \Phi_{B}(x-x_j(t)) R_j^{p-1}\Big) \e^2 
\\ &&
+\sum_{j=1}^N \int \Phi_{B}(x-x_j(t))
 (c(t,x)-c_j(t)) \e^2
\\&&
+\int \Big(1-\sum_{j=1}^N \Phi_{B}(x-x_j(t))\Big)
( \e^2_x  + c(t,x) \e^2 ).
\eee
Next, we make the following observations:

(i) By (\ref{singlelocal}), we have
$\forall j$,
$$\int \Phi_{B}(x-x_j(t))
\left(\e^2_x - p R_j^{p-1} \e^2 + c_j(t)\e^2\right) 
\ge {\l_1\over 4} \int  \Phi_{B}(x-x_j(t))
(\e^2_x + c_j(t)\e^2).
$$

(ii) Since $\Phi_B(x)=1$ for $|x|<B$, by the decay properties
of $Q$, we have
\bee
0\le R^{p-1}-\sum_{j=1}^N \Phi_{B}(x-x_j(t)) R_j^{p-1}
\le |R|_{L^\infty(|x-x_j(t)|>B)}^{p-1} + C \sum_{j\not =k} R_j R_k
\le C e^{-\sqrt{\sigma_0} B}.
\eee

(iii) Note that 
$c(t,x)=\sum_{j=1}^N c_j(t) \varphi_j(t,x),$
where $\varphi_1(t,x)=1-\psi(x-m_2(t)),$
for $j\in\{2,\ldots,N-1\}$,
$\varphi_j(t,x)=\psi(x-m_{j}(t))-\psi(x-m_{j+1}(t))$ and
$\varphi_N(t,x)=\psi(x-m_N(t)).$
Since $\Phi_B(x)\le 3 e^{-{|x|\over B}}$, 
by the properties of $\psi$, and 
$|m_j(t)-x_j(t)|\ge L_4/2\ge 2 k B,$  we obtain
\bee
\left| \Phi_B(x-x_j(t)) (c(t,x)-c_j(t))\right|
&\le& |c(t,x)-c(t)|_{L^\infty(|x-x_j(t)|\le k B)}
+C e^{- k }
\\ &\le& C e^{-\sqrt{\sigma_0} k B/2} +C e^{- k} .
\eee

(iv) $1-\sum_{j=1}^N \Phi_{B}(x-x_j(t))\ge 0$.

Therefore, with $\l_0={1\over 2}
 \min({\l_1\over 4},{\l_1\over 4}\sigma_0,1,\sigma_0)$,
for $B$ and $k$ large enough,
\bee
\int \e^2_x  - p R^{p-1} \e^2 + c(t,x) \e^2 &\ge &
2 \l_0 \int (\e_x^2+\e^2) -C \left(e^{-\sqrt{\sigma_0} B/2} + e^{- k}
\right) \int \e^2
\\
&\ge & \l_0 \int (\e_x^2+\e^2).
\eee
Thus  the proof of Lemma \ref{QUAD} is complete.

\section{Proof of the stability in the energy space}

This section is devoted to the proof of stability result. The proof is by 
a priori estimate.

Let $0<c_1^0<\ldots<c_N^0$,
$\sigma_0={1\over 2}\min(c_1^0,c_2^0-c_1^0,c_3^0-c_2^0,\ldots,c_N^0-c_{N-1}^0)
$ and  $\gamma_0=\sqrt{\sigma_0}/16$.
For $A_0,L,\alpha>0$, we define 
\be\label{defofv}
{\cal V}_{A_0}(L,\alpha)=\Big\{ u\in H^1(\R);
\inf_{x_j -x_{j-1} \ge L} 
\Big|u-\sum_{j=1}^N Q_{c_j^0} (.-x_j)\Big|_{H^1} \le A_0 \Big( \alpha
+e^{-\gamma_0 L/2} \Big) \Big\}.
\ee

We want to prove that there exists  
$A_0>0$, $L_0>0$, and $\alpha_0>0$ such that,
$\forall u_0\in H^1(\R)$, if for some $L>L_0$, $\alpha<\alpha_0$,
$
\left| u_0-\sum_{j=1}^N Q_{c_j^0}(.-x_j^0) \right|_{H^1}\le\alpha,
$
where $x_j^0>x_{j-1}^0+L$, then
$
\forall t\ge 0$, $
u(t)\in {\cal V}_{A_0}(L,\alpha)
$
(this  proves the stability result in $H^1$).
By a standard continuity argument, 
it is a direct consequence of the following proposition.

\begin{proposition}[A priori estimate]\label{MAINPROP}
There exists  
$A_0>0$, $L_0>0$, and $\alpha_0>0$ such that,
for all $u_0\in H^1(\R)$, if 
\be\label{initprop} 
\Big| u_0-\sum_{j=1}^N Q_{c_j^0}(.-x_j^0) \Big|_{H^1}\le\alpha,
\ee
where $L>L_0$, $0<\alpha<\alpha_0$, $x_j^0>x_{j-1}^0+L$, 
and if for $t^*>0$,
\be\label{aprioriprop}
\forall t\in [0,t^*],\quad 
u(t)\in {\cal V}_{A_0}(L,\alpha),
\ee
then
\be\label{conclprop}
\forall t\in [0,t^*],\quad 
u(t)\in {\cal V}_{A_0/2}(L,\alpha).
\ee
\end{proposition}

Note that $A_0$, $L_0$ and $\alpha>0$ are independent
of $t^*$.

{\bf Proof of Proposition \ref{MAINPROP}.}\quad 
Let $A_0>0$ to be fixed later.
First, for $0<\alpha_0<\alpha_I(A_0)$ and $L_0>L_I(A_0)>L_1$, we have
\be\label{allsmall0}
A_0\left(\alpha_0+ e^{-\gamma_0 L_0/2}\right)\le \alpha_1,
\ee
where $\alpha_1$ and $L_1$ are   defined in Lemma \ref{DECOMPOSITION}.
Therefore, by (\ref{aprioriprop}) and Lemma \ref{DECOMPOSITION}, 
there exist
$c_j:[0,t^*]\goto (0,+\infty)$,
$x_j:[0,t^*]\goto \R$, such that
\be\label{proofdefofeps}
\e(t,x)=u(t,x)-\sum_{j=1}^N R_j(t,x),
\quad \hbox{where} \quad R_{j}(t,x)=Q_{c_j(t)}(x-x_j(t)),
\ee
satisfies $\forall j$, $\forall t\in [0,t^*],$
\be\label{proofortho}
\int R_j(t) \e(t)=
\int (R_j(t))_x \e(t)=0,
\ee
\be\label{smallapriori}
|c_j(t)-c_j^0|+|\dot c_j|+|\dot x_j-c_j^0|+
|\e(t)|_{H^1} \le K_1 (A_0+1)  \left(\alpha_0 +e^{-\gamma_0 L_0}\right).
\ee
Note that by \eqref{initprop}, Lemma \ref{MODULATION} (see Appendix)
and assumptions of the proposition,
\be\label{important}
|\e(0)|_{H^1}+\sum_{j=1}^N|c_j(0)-c_j^0|\le 
K_1 \alpha,\quad x_j(0)-x_{j-1}(0)\ge {L\over 2}.
\ee
From \eqref{smallapriori} and (\ref{important}),
for  $\alpha_0<\alpha_{II}(A_0)$
 and $L_0>L_{II}(A_0)>2 \max(L_2,L_3,L_4)$
($L_2$, $L_3$ and $L_4$ are defined in Lemmas \ref{th:2} and \ref{QUAD}), 
we have $\forall t\in [0,t^*],$
\bea
&& c_1(t)\ge \sigma_0, \quad \dot x_1(t) \ge \sigma_0,
\quad c_j(t)-c_{j-1}(t)\ge \sigma_0, \quad \dot x_j(t)-\dot x_{j-1}(t)
\ge \sigma_0, \label{allsmall2}\\
&& x_j(t)-x_{j-1}(t)\ge L/2\ge \max(L_3,L_4),
\quad 
|\e(t)|_{H^1}\le {1\over 2} \left({\sigma_0\over 8}\right)^{1\over p-1}.
\label{allsmall3}
\eea
Therefore, we can apply Lemmas \ref{th:3}, \ref{th:2} and \ref{QUAD} 
for all $t\in [0,t^*]$.

Let $\alpha_0=\min(\alpha_I(A_0),\alpha_{II}(A_0))$ and
$L_0=\max(L_{I}(A_0),L_{II}(A_0)).$
Now, our objective is to give a uniform upper bound on $|\e(t)|_{H^1}$
and $|c_j(t)-c_j(0)|$ on $[0,t^*]$ improving \eqref{smallapriori}
for $A_0$ large enough.

\medskip

In the next lemma, 
we first obtain a control of  the variation of  $c_j(t)$ 
which is quadratic in $|\e(t)|_{H^1}$.
This is the key step of the stability result, based on monotonicity 
property of the local $L^2$ norm and energy constraints.
It is  essential at this point to have chosen by the modulation
$\int  R_j \e=0$. 

\begin{lemma}[Quadratic control of the variation of $c_j(t)$]\label{PROOFL1}
There exists $K_4>0$ 
independent of $A_0$, such that, $\forall t\in [0,t^*]$, 
\be\label{proofl1}
\sum_{j=1}^N \left|c_j(t) - c_j(0) \right|
\le K_4 \left( |\e(t)|_{H^1}^2 +|\e(0)|_{H^1}^2 +e^{-\gamma_0 L}\right).
\ee
\end{lemma}

{\bf Proof.}

{\sl Step 1.} Energetic control.
\quad Let $\beta={2\over p-1}$. There exists $C>0$ such that
\bea
\left| \sum_{j=1}^N  c_j(0) \bkt{c_j^{\beta-1/2}(t)
- c_j^{\beta-1/2}(0) } \right|
&\le& C\left( |\e(t)|_{H^1}^2 +|\e(0)|_{H^1}^2 +e^{-\gamma_0 L}\right)
 \nonumber \\ &&\quad 
+C \sum_{j=1}^N \bkt{c_j(t)
- c_j(0) }^2. \label{encorebis}
\eea

Let us prove (\ref{encorebis}).
By \eqref{29}, we have
\be\label{29bisbis}
\left| \sum_{j=1}^N \bkt{E(R_j(t))-E(R_j(0))}\right|
\le C\left( |\e(t)|_{H^1}^2 +|\e(0)|_{H^1}^2 +e^{-\gamma_0 L}\right).
\ee
Since
$ E(Q_c)=- {\k\over 2} \,  c^{\beta+1/2}  \int Q^2,$
where $\k=\frac {5-p}{p+3}$,
we have
$$
- \sum_{j=1}^N \bkt{E(R_j(t))-E(R_j(0))}
= {\k\over 2}  \left(\int Q^2\right)  \sum_{j=1}^N \bkt{c_j^{\beta+1/2}(t)
- c_j^{\beta+1/2}(0)}.$$
By linearization, we have
$c_j^{\beta+1/2}(t) - c_j^{\beta+1/2}(0)=
{2\beta+1\over 2\beta-1} c_j(0) \bkt{c_j^{\beta-1/2}(t)
- c_j^{\beta-1/2}(0)}
+O\Big(\bkt{c_j(t)
- c_j(0) }^2\Big).$ Note that ${2\beta+1\over 2\beta-1}={1\over \k}$.
Therefore,
\bea
& & \left| \sum_{j=1}^N \bkt{E(R_j(t))-E(R_j(0))}
+{1\over 2}  \left(\int Q^2\right)
 \sum_{j=1}^N  c_j(0) \bkt{c_j^{\beta-1/2}(t)
- c_j^{\beta-1/2}(0) } \right|  \nonumber \\ &&
\qquad \le  C \sum_{j=1}^N \bkt{c_j(t)
- c_j(0) }^2,
\label{34}
\eea
and from (\ref{29bisbis}), we obtain (\ref{encorebis}).

\medskip

{\sl Step 2.} $L^2$ mass monotonicity at the right of every soliton.
Let 
$$ d_j(t)= \sum_{k=j}^N  c_k^{\beta-1/2}(t).$$
We claim
\be\label{proof3}  \left(\int Q^2\right)  \left|d_j(t)-d_j(0)\right|
 \le - \left(\int Q^2\right)  (d_j(t)-d_j(0)) + C \Big[ 
\int \e^2(0) + e^{-\gamma_0 L}\Big].
\ee

Let us prove (\ref{proof3}).
Recall that using the notation of section \S 2.3, we have
$${\cal I}_j(t)\le {\cal I}_j(0) + K_3 e^{-\gamma_0 L},\quad
\hbox{where} \quad {\cal I}_j(t)=\int \psi(x-m_j(t)) u^2(t,x) dx.$$
Since $\int R_j^2(t)=c_j^{\beta-1/2} (t) \int Q^2$, $\int R_j(t) \e(t)=0$,
by similar calculations as in Lemma \ref{th:3}, we have
\be\label{proof1}
\left| {\cal I}_j(t) -\left(\int Q^2\right) d_j(t)
 - \int \psi(.-m_j(t)) \e^2(t)
\right|
\le C e^{-\gamma_0 L}.
\ee
Therefore,
\be\label{proof2}
\left(\int Q^2\right) (d_j(t)-d_j(0)) \le \int \psi(.-m_j(0)) \e^2(0)-\int \psi(.-m_j(t)) \e^2(t)+ C e^{-\gamma_0 L}.
\ee
Since the second term on the right hand side is negative,
(\ref{proof3}) follows easily.
Note that by conservation of the $L^2$ norm $\int u^2(t)=\int u^2(0)$
and 
$$\int u^2(t)=\int R^2(t)+\int \e^2(t)+2\int R(t)\e(t)
=\int R^2(t)+\int \e^2(t)= d_1(t)+\int \e^2(t)+O(e^{-\gamma_0 L}),
$$ 
we obtain
\be\label{proof4} \left(\int Q^2\right) (
d_1(t)-d_1(0)) \le \int  \e^2(0) - \int \e^2(t) + C e^{-\gamma_0 L}.
\ee

\medskip

{\sl Step 3.} Resummation argument.
\quad 
By Abel transform, we have
\bea &&
\sum_{j=1}^N  c_j(0) \bkt{c_j^{\beta-1/2}(t)
- c_j^{\beta-1/2}(0) } \nonumber\\ && \qquad = \sum_{j=1}^{N-1} c_j(0) \bkt{ d_j(t)-d_{j+1}(t)
-(d_j(0)-d_{j+1}(0))} 
+c_N(0) \bkt{d_N(t)-d_N(0)} \nonumber \\ &&\qquad 
= c_1(0) \bkt{d_1(t)-d_1(0)}
  +\sum_{j=2}^N  (c_j(0)-c_{j-1}(0)) (d_j(t)-d_j(0)).
\label{proof5} \eea
Therefore, by step 1, 
\bea  &&
-\left(c_1(0) \bkt{d_1(t)-d_1(0)}
  +\sum_{j=2}^N  (c_j(0)-c_{j-1}(0)) (d_j(t)-d_j(0)) \right) \nonumber \\ &&\qquad 
\le C\left( |\e(t)|_{H^1}^2 +|\e(0)|_{H^1}^2 +e^{-\gamma_0 L}\right) + 
 C \sum_{j=1}^N \bkt{c_j(t) - c_j(0) }^2.\label{unnom}\eea
Since $c_1(0)\ge\sigma_0$, 
$c_j(0)-c_{j-1}(0)\ge \sigma_0$, by \eqref{proof3}, we have
\bee \sigma_0 \sum_{j=1}^N \left|d_j(t) - d_j(0) \right|
&\le&  c_1(0) |d_1(t)-d_1(0)|
  +\sum_{j=2}^N  (c_j(0)-c_{j-1}(0)) |d_j(t)-d_j(0)| 
\\
&\le&  - \left[ c_1(0) \bkt{d_1(t)-d_1(0)}
  +\sum_{j=2}^N  (c_j(0)-c_{j-1}(0)) (d_j(t)-d_j(0)) \right]
\\ &&  + C\int \e^2(0)+
Ce^{-\gamma_0 L}.
\eee
Thus, by (\ref{unnom}), we have
\bee \sum_{j=1}^N \left|d_j(t) - d_j(0) \right|\le C\left( |\e(t)|_{H^1}^2 +
|\e(0)|_{H^1}^2 +e^{-\gamma_0 L}\right) + 
 C \sum_{j=1}^N \bkt{c_j(t) - c_j(0) }^2.
\eee

Since
\bee |c_j(t)-c_j(0)|&\le& C|c_j^{\beta-1/2}(t) -c_j^{\beta-1/2}(0)| \\
&\le& C(|d_j(t)-d_j(0)|+|d_{j+1}(t)-d_{j+1}(0)|),
\eee
we obtain, 
\bee \sum_{j=1}^N \left|c_j(t) - c_j(0) \right| 
\le 
C\left( |\e(t)|_{H^1}^2 +|\e(0)|_{H^1}^2 +e^{-\gamma_0 L}\right)
 +C \sum_{j=1}^N \bkt{c_j(t) - c_j(0) }^2.
\eee
Choosing a smaller $\alpha_0(A_0)$ and a larger
$L_0(A_0)$, by (\ref{smallapriori}), we assume 
$C  |c_j(t) - c_j(0)| \le 1/2$ and so
\be\label{proof6}
\sum_{j=1}^N
\left|c_j(t) - c_j(0) \right|
\le C\left( |\e(t)|_{H^1}^2 +|\e(0)|_{H^1}^2 +e^{-\gamma_0 L}\right).
\ee
Thus, Lemma \ref{PROOFL1} is proved.

\medskip

Now, we prove the following lemma, giving uniform control on $|\e(t)|_{H^1}$
on $[0,t^*]$. 

\begin{lemma}[Control of $|\e(t)|_{H^1}$]\label{proofl2}
There exists 
$K_5>0$ independent of $A_0$, such that, $\forall t\in [0,t^*],$
$$|\e(t)|_{H^1}^2
\le K_5 \left( |\e(0)|^2_{H^1}+e^{-\gamma_0 L}\right).$$
\end{lemma}

{\bf Proof.} \quad It follows from direct calculation on the energy, and 
the previous estimates obtained by  Abel transform, freezing the $c_j(t)$ at
the first order.

By (\ref{29}), (\ref{34}), (\ref{proof5}) and (\ref{proof6}), we have
\bee
&&\frac 12 \int \e_x^2(t) - p  R^{p-1}(t) \e^2(t)
 \\ &&\quad \le 
-\sum_{j=1}^N \bkt{E(R_j(t))-E(R_j(0))}
 + K_2\left( |\e(0)|^2_{H^1}+|\e(t)|_{H^1}^3
+e^{-\gamma_0 L}\right) \\
&& \quad \le 
{1\over 2}  \left(\int Q^2\right) \sum_{j=1}^N  c_j(0) \bkt{c_j^{\beta-1/2}(t)
- c_j^{\beta-1/2}(0) } + C \sum_{j=1}^N \bkt{c_j(t)
- c_j(0) }^2 \\ && \quad \quad +K_2\left( |\e(0)|^2_{H^1}+|\e(t)|_{H^1}^3
+e^{-\gamma_0 L}\right) \\
&& \quad \le  {1\over 2} \left(\int Q^2\right) \bkt{
 c_1(0) \bkt{d_1(t)-d_1(0)}
  +\sum_{j=2}^N  (c_j(0)-c_{j-1}(0)) (d_j(t)-d_j(0))}\\ &&\quad \quad
+C \left( |\e(0)|^2_{H^1}+|\e(t)|_{H^1}^3
+e^{-\gamma_0 L}\right).
\eee
Therefore,  using \eqref{proof2} and \eqref{proof4}, and again
Lemma \ref{PROOFL1}, we have
\bea
\int \e_x^2(t) - p  R^{p-1}(t) \e^2(t) 
& \le &
- \Big( c_1(0) \int \e^2(t) +\sum_{j=2}^N 
(c_j(0)-c_{j-1}(0)) \int \psi(x-m_j(t)) \e^2(t)\Big)\nonumber
\\ && + C \left( |\e(0)|^2_{H^1}+|\e(t)|_{H^1}^3
+e^{-\gamma_0 L}\right)\nonumber 
\\
&\le &
- \int c(t,x) \e^2(t)+ C \left( |\e(0)|^2_{H^1}+|\e(t)|_{H^1}^3
+e^{-\gamma_0 L}\right)\label{dernier}
\eea
where
$c(t,x)=  c_1(t) + \sum_{j=2}^N ( c_j(t)- c_{j-1}(t))\psi(x-m_j(t)).$

By Lemma \ref{QUAD},
$$\int \e_x^2(t) - p  R^{p-1}(t) \e^2(t) 
+c(t,x) \e^2(t) \ge \l_0 |\e(t)|_{H^1}^2.$$
Therefore, from (\ref{dernier}), we obtain
$$|\e(t)|_{H^1}^2 \le C \left( |\e(0)|_{H^1}^2+|\e(t)|_{H^1}^3
+e^{-\gamma_0 L}\right),$$
and so
$$|\e(t)|_{H^1}^2 \le K_5 \left(|\e(0)|^2_{H^1}+e^{-\gamma_0 L}\right),$$
for some constant $K_5>0$, independent of $A_0$.
Thus Lemma \ref{proofl2} is proved.

\bigskip

We conclude the proof of proposition \ref{MAINPROP} and
of the stability result.
By (\ref{important}) and Lemmas \ref{PROOFL1} and \ref{proofl2}, we have
\bee  &&
\Big|u(t)-\sum_{j=1}^N Q_{c_j^0}(x-x_j(t))\Big|_{H^1} \\ &&\qquad \le
\Big|u(t)-\sum_{j=1}^N R_j(t)\Big|_{H^1} 
+\Big|\sum_{j=1}^N R_j(t)-\sum_{j=1}^N  Q_{c_j^0}(x-x_j(t))\Big|_{H^1}\\
&&\qquad \le  |\e(t)|_{H^1} + C\sum_{j=1}^N|c_j(t)-c_j^0| \\
&& \qquad \le 
 |\e(t)|_{H^1} + C\sum_{j=1}^N|c_j(t)-c_j(0)|+ C\sum_{j=1}^N|c_j(0)-c_j^0|\\
&& \qquad \le  |\e(t)|_{H^1} + C K_4 (|\e(0)|^2_{H^1}+e^{-\gamma_0 L})
 + CK_1 \alpha \\
&&\qquad \le   K_6 \left( \alpha +e^{-\gamma_0 L/2}\right),\eee
where $K_6>0$ is a constant independent of $A_0$.

Choosing  $A_0=4 K_6$, we complete the proof of Proposition \ref{MAINPROP}
and thus the proof of Theorem \ref{MAINTH} (i).

\section{Proof of the asymptotic stability result}
This section is devoted to the proof of the asymptotic stability
result (Theorem \ref{MAINTH} (ii)).

\subsection{Asymptotic stability around the solitons}
In this subsection, we prove the following asymptotic result on $\e(t)$
as $t\goto +\infty$.

\begin{proposition}[Convergence around solitons, $p=2,3,4$]\label{CONVAROUND}
Under the assumptions of Theorem \ref{MAINTH}, the following
is true:\\ 
(i) Convergence of $\e(t)$: \quad $\forall j\in \{1,\ldots,N\}$,
\be\label{convaround}
\e(t,.+x_j(t))\rightharpoonup 0 \quad \hbox{in $H^1(\R)$ as
$t\goto +\infty$}.
\ee 
\noindent (ii) Convergence of geometric parameters: \quad  there exists
$0<c_1^{+\infty}<\ldots<c_N^{+\infty}$, such that
$$c_j(t)\goto c_j^{+\infty},\quad 
\dot x_j(t) \goto c_j^{+\infty} \quad  
\hbox{as $t\goto +\infty$}.$$
\end{proposition}

The proof of this result is very similar to the proof of the asymptotic
stability of a single soliton in Martel and Merle \cite{MM3}
for the subcritical case (see also the previous paper \cite{MM2}
concerning the critical case $p=5$). The proof is based on the 
following rigidity result of solutions of \eqref{kdvp} around solitons.

\medskip

\noindent{\bf Theorem} ({\bf Liouville property close to $R_{c_0}$ for $p=2,3,4$
\cite{MM3}})
{\it
Let $p=2,3 $ or $4$, and let $c_0>0$.
 Let $u_0\in H^1(\R)$, and let $u(t)$ be the solution of (\ref{kdvp}) 
for all time $t\in \R$.
There exists $\alpha_0>0$ such that
if $|u_0-R_{c_0}|_{H^1}<\alpha_0$, and if
there exists $y(t)$ such that
\begin{equation}\label{l2compact}\forall \delta_0>0,\exists A_0>0/
\forall t\in \R,
\quad \int_{|x|>A_0} u^2(t,x+y(t))
dx\le \delta_0, \quad \hbox{($L^2$ compactness)}
\end{equation} 
then there exists $c^*>0,$ $ x^*\in  \R$ such that 
$$\forall t\in \R, \forall x\in \R,\quad u(t,x)=Q_{c^*}(x-x^*-c^* t).$$
}

{\bf Proof of Proposition \ref{CONVAROUND} (i).} \quad 
Consider a solution $u(t)$ satisfying the assumptions of Theorem \ref{MAINTH}.
Then, by \S 3, we known that $u(t)$ is uniformly close in $H^1(\R)$ 
to the superposition of $N$ solitons 
for all time $t\ge 0$.
With the decomposition introduced in section \S 2, it is equivalent that
$\e(t)$ is uniformly small in $H^1(\R)$ and $\sum_{j=1}^N |c_j(t)-c_j(0)|$
is uniformly small. Therefore, we can assume that, $\forall t\ge 0$,
$$c_1(t)\ge \sigma_0,\quad c_j(t)-c_{j-1}(t)\ge \sigma_0.$$

\medskip

The proof of Proposition \ref{CONVAROUND} is by contradiction.
Let $j\in \{1,\ldots,N\}$.
Assume that for some sequence $t_n\goto +\infty$, we have
$$\e(t_n,.+x_j(t_n))\not \rightharpoonup 0 \quad \hbox{in $H^1(\R)$ as
$t\goto +\infty$}.
$$
Since $0<\sigma_0<c_j(t)<\overline c$ 
and $|\e(t)|_{H^1}\le C$ for all $t\ge 0$, there exists
$\et_0\in H^1(\R)$, $\et_0\not \equiv 0$,
and $\ct_0>0$ such that for a subsequence of $(t_n)$,
still denoted $(t_n)$, we have
\be\label{contra}
\e(t_n,.+x_j(t_n)) \rightharpoonup \et_0 \quad \hbox{in $H^1(\R)$,} 
\quad c_j(t_n)\goto \ct_0
\quad \hbox{as $n\goto +\infty$}.
\ee
Moreover, by weak convergence and  stability result, 
$|\et_0|_{H^1}\le \sup_{t\ge 0} |\e(t)|_{H^1}\le C (\alpha_0
+e^{-\gamma_0 L_0})$, and therefore $|\et_0|_{H^1}$
is as small as we want by taking $\alpha_0$ small and $L_0$ large.

\medskip

Let now $\ut(0)=Q_{\ct_0}+\et_0$, and let $\ut(t)$ be the global
solution of \eqref{kdvp} for $t\in \R$, with $\ut(0)$ as initial data.
Let $\xt(t)$ and $\ct(t)$ be the geometrical parameters associated to
the solution $\ut(t)$ (apply the modulation theory for a solution
close to a single soliton).

We claim that the solution $\ut(t)$ is $L^2$ compact in the sense
of \eqref{l2compact}.

\begin{lemma}[$L^2$ compactness of the asymptotic solution]
\label{COMPACTBIS}
\begin{equation}\label{l2compactbis}
\forall \delta_0>0,\exists A_0>0/\forall t\in \R,
\quad \int_{|x|>A_0} \ut^2(t,x+\xt(t))
dx\le \delta_0.
\end{equation} 
\end{lemma}

Assuming this lemma, we finish the proof of Proposition 
\ref{CONVAROUND} (i).
Indeed, by choosing $\alpha_0$ small enough and $L_0$ large enough,
we can apply the Liouville theorem to $\ut(t)$. Therefore, there
exists $c^*>0$ and $x^*\in \R$, such that
$\ut(t)=Q_{c^*}(x-x^*-c^* t).$
In particular, $\ut(0)= Q_{\ct_0}+\et_0= Q_{c^*}(x-x^*).$
Since by weak convergence
$\int \et_0 (Q_{c_0})_x=0$,
we have easily $x^*=0$. Next, since 
$ \int \et_0 Q=0$, we have $c^*=\ct_0$ and so 
$\et_0\equiv 0$.
This is a contradiction.

Thus Proposition \ref{CONVAROUND} (i) is proved assuming
Lemma \ref{COMPACTBIS}.
The proof of Lemma \ref{COMPACTBIS} 
is based only on arguments of monotonicity of the
$L^2$ mass in the spirit of \cite{MM3}, \cite{MM4}.

\medskip

{\bf Proof of Lemma \ref{COMPACTBIS}.}
\quad 
We use the function $\psi$ introduced in \S 2.2.
For $y_0>0$, we introduce two quantities:
\be\label{LR}
J_L(t)=\int (1-\psi(x-(x_j(t)-y_0))) u^2(t,x) dx,\quad
J_R(t)=\int \psi(x-(x_j(t)+y_0)) u^2(t,x) dx.
\ee

The strategy of the proof is the following.
We prove first that $J_L(t)$ is almost increasing
and $J_R(t)$ is almost decreasing in time.
Then, assuming  by contradiction that $\ut(t)$ is not
$L^2$ compact, using the convergence of $u(t)$ to 
$\ut(t)$ for all time, we prove  that the $L^2$ norm of
$u(t)$ in the compact set $[-y_0,y_0]$,
for $y_0$ large enough, oscillates between two 
different  values. This proves that there are infinitely many  
transfers of mass
from the right hand side of the soliton $j$ to the 
left hand side of the soliton $j$. This is of course impossible since
the $L^2$ norm of $u(t)$ is finite.

\medskip

{\bf Step 1.}\quad Monotonicity on the right and on the left of
a soliton. \quad 
We claim

\medskip

{\sl Claim.}\quad There exists $C_1,y_1>0$
such that 
$\forall y_0>y_1$, $\forall t'\in [0,t]$,
\be\label{monotonie}
J_L(t)\ge J_L(t')-C_1 e^{-\gamma_0 y_0},\quad 
J_R(t)\le J_R(t')+C_1 e^{-\gamma_0 y_0}.
\ee

\medskip

We prove this claim.\quad First note that it is sufficient to prove
(\ref{monotonie}) for  $J_L(t)$. Indeed,
since $u(-t,-x)$ is also solution of (\ref{kdvp}),
and since $1-\psi(-x)=\psi(x)$, 
we can argue backwards in time (from $t$ to $t'$) to obtain the
result for $J_R(t)$. By using the
same argument as in Lemma \ref{th:2}, we prove easily,
for $y_0$ large enough, for all $0<t'<t$,
\bee \int \psi(.-(x_j(t)-y_0-{\textstyle{\sigma_0\over 2}}(t-t'))) u^2(t) 
&\le & \int \psi(.-(x_j(t')-y_0)) u^2(t') + C_1 e^{-\gamma_0 y_0} \\
&\le & \int u^2(t') - J_L(t') + C_1 e^{-\gamma_0 y_0}.
\eee
Since $\int u^2(t)=\int u^2(t')$ and 
$$\int u^2(t)-J_L(t)=\int \psi(.-(x_j(t)-y_0)) u^2(t)
\le \int \psi(.-(x_j(t)-y_0-{\textstyle {\sigma_0\over 2}}(t-t'))) u^2(t),$$
we obtain the result.

\medskip

{\bf Step 2.}\quad Conclusion of the proof. 
\quad 
Recall from \cite{MM3} that we have stability of (\ref{kdvp})
by weak convergence in $H^1(\R)$ in the following sense
\be\label{stabconv}
\forall t\in \R,\quad
u(t+t_n,.+x_j(t+t_n))\longrightarrow  \ut(t,.+\xt(t))\quad 
\hbox{in $L^2_{loc}(\R)$ as $n\goto +\infty$.}
\ee
This was proved in \cite{MM3} by using the fact
that the Cauchy problem for (\ref{kdvp}) is well posed
both in $H^1(\R)$ and in $H^{s^*}(\R)$, for some
$0<s^*<1$, for any $p=2,3,4$ (see \cite{KPV}).

\medskip

We prove Lemma \ref{COMPACTBIS} by contradiction.
Let
$$m_0=\int \ut^2(0)=\int  \ut^2(t).$$
Assume that there exists
$\delta_0>0$ such that for any $y_0>0$, there exists $t_0(y_0)\in \R$,
such that
\be\label{loss}\int_{|x|< 2 y_0} \ut^2(t_0(y_0),x+\xt(t_0(y_0))) dx 
\le m_0 -\delta_0.\ee
Fix $y_0>0$ large enough so that
\be\label{loss3}
\int (\psi(x+y_0)-\psi(x-y_0))\ut^2(0,x) dx
\ge m_0 -{1\over 10} \delta_0,
\ee
$$
C_1 e^{-\gamma_0 y_0}+m_0 \sup_{|x|>2 y_0}\{\psi(x+y_0)-\psi(x-y_0)\} 
\le {1\over 10} \delta_0.
$$
Assume that $t_0=t_0(y_0)>0$ and, by possibly considering 
a subsequence of $(t_n)$, that $\forall n$, $t_{n+1}\ge t_n+t_0$.

Observe that, since $0<\psi<1$ and 
$\psi'>0$,  by the choice of $y_0$ and (\ref{loss}), we have 
\bea &&
\int (\psi(x-(\xt(t_0)-y_0))-\psi(x-(\xt(t_0)+y_0)))
\ut^2(t_0,x)dx 
 \nonumber \\&& \quad \le 
\int_{|x|< 2 y_0} \ut^2(t_0,x+\xt(t_0)) dx 
+m_0 \sup_{|x|>2 y_0}\{\psi(x+y_0)-\psi(x-y_0)\} 
\nonumber \\&& \quad \le 
\int_{|x|< 2 y_0} \ut^2(t_0,x+\xt(t_0)) dx 
+{1\over 10} {\delta_0}\le m_0 -{9\over 10} \delta_0.
\label{loss2}
\eea
Then, by (\ref{loss3}), (\ref{loss2})
and  (\ref{stabconv}),
there exists $N_0>0$ large enough
so that $\forall n\ge N_0,$
\be\label{osc2}
\int (\psi(x-(x_j(t_n)-y_0))-\psi(x-(x_j(t_n)+y_0))) u^2(t_n,x) dx 
\ge m_0 -{1\over 5} \delta_0.
\ee
\be\label{osc1}
\int (\psi(x-(x_j(t_n+t_0)-y_0))-\psi(x-(x_j(t_n+t_0)+y_0)))
  u^2(t_n+t_0,x) dx  \le m_0 -{4\over 5} \delta_0.\ee

Recall that from Step 1, and the choice of $y_0$, we have
$
J_R(t_n+t_0)\le J_R(t_n) +{1\over 10}\delta_0.$
Therefore, by conservation of the $L^2$ norm and (\ref{osc1}),
 (\ref{osc2}), we have
$$J_L(t_n+t_0)\ge J_L(t_n) +{1\over 2} \delta_0.$$
Since $J_L(t_{n+1})\ge J_L(t_n+t_0)-{1\over 10} \delta_0$ by Step 1,
we finally obtain
$$\forall n\ge N_0,\quad
J_L(t_{n+1})\ge J_L(t_n) + {2\over 5} \delta_0.$$
Of course, this is a contradiction.
Thus the proof of Lemma \ref{COMPACTBIS} is complete.

\medskip

{\bf Proof of Proposition \ref{CONVAROUND} (ii).}
The proof is similar to the proof of Proposition 3 in \cite{MM3}.
It follows again from monotonicity arguments and the fact that 
we consider the subcritical case $1<p<5$.

Let $\delta>0$ be arbitrary.
Since 
$\int R_j^2(t)=c_j^{5-p\over 2(p-1)}(t) \int Q^2$ and 
$\e(t,.+x_j(t))\goto  0$ in $L^2_{loc}$ as $t\goto +\infty$, 
there exists $T_1(\delta)>0$ and $y_1(\delta)$ such that
$\forall t>T_1(\delta)$, $\forall y_0>y_1(\delta)$,
$$\left| \int (\psi(x-(x_j(t)-y_0))-\psi(x-(x_j(t)+y_0))) u^2(t,x) dx
-c_j^{5-p\over 2(p-1)}(t) \int Q^2 \right|
\le \delta.$$

By Step 1 of  the proof of Lemma \ref{COMPACTBIS}, there exists
$y_2(\delta)$, such that 
we have, for all $0<t'<t$, $\forall y_0>y_2(\delta)$,
$$J_L(t)\ge J_L(t') - \delta,\quad J_R(t)\le J_R(t')+\delta.
$$
Fix $y_0=\max(y_1(\delta),y_2(\delta))$, it follows that 
there exists  $T_2(\delta)$, $J_L^{+\infty}\ge 0$ and
$J_R^{+\infty}\ge 0$ such that
$$\forall t\ge T_2(\delta),
\quad
|J_L(t)-J_L^{+\infty}|\le 2 \delta,\quad 
|J_R(t)-J_R^{+\infty}|\le 2 \delta.
$$

Therefore, by conservation of $L^2$ mass, we have,
for all $0<\max(T_1,T_2)<t'<t$,
$$\left|c_j^{5-p\over 2(p-1)}(t)- c_j^{5-p\over 2(p-1)}(t')\right|
\le C\, \delta.$$  
Since $\delta$ is arbitrary, it follows that 
$c_j^{5-p\over 2(p-1)}(t)$ 
has a limit as $t\goto +\infty$.
Thus there exists $c_j^{+\infty}>0$ such that
$c_j(t)\goto c_j^{+\infty}$ as $t\goto +\infty$.
The fact that $\dot x_j(t)\goto c_j^{+\infty}$ is a direct
consequence of (\ref{faux13}).

\subsection{Asymptotic behavior on $x>c t$}

In this subsection, using the same argument of monotonicity of
$L^2$ mass, we prove the following proposition.

\begin{proposition}[Convergence for $x>c_1^0 t/10$]\label{ALLX}
Under the assumptions of Theorem \ref{MAINTH}, the following is
true
\be\label{allx}
|\e(t)|_{L^2(x>c_1^0 t /10)} \goto 0\quad \hbox{as $t\goto +\infty$}.
\ee
\end{proposition}

{\bf Proof.}\quad 
By arguing backwards in time (from $t$ to 0) and using the conservation
of $L^2$ norm, we have
$$
\int \psi(.-(x_N(t)+y_0)) u^2(t) 
\le \int \psi(.-(x_N(0)+{\textstyle {\sigma_0\over 2}} 
t +y_0)) u^2(0) + C_1 e^{-\gamma_0 y_0}.
$$
Therefore, 
$$\int_{x>x_N(t)+y_0} \e^2(t) \le
 2 \int \psi(.-(x_N(0)+{\textstyle {\sigma_0\over 2}}
 t +y_0)) u^2(0)
 + C e^{-\gamma_0 y_0}.
$$
Since for fixed $y_0$, $\int_{x_N(t)<x<x_N(t)+y_0} \e^2(t)
\goto 0$ as $t\goto +\infty$, we conclude
$\int_{x>x_N(t)} \e^2(t) \goto 0$  as $t\goto +\infty.$

Now, let us prove $\int_{x>x_j(t)} \e^2(t) \goto 0$ as $t\goto +\infty$
by backwards induction on $j$.
Assume that for $j_0\in \{2,\ldots,N\}$, we have
$\int_{x>x_{j_0}(t)} \e^2(t) \goto 0$ as $t\goto +\infty$.
For $t\ge 0$ large enough, there exists $0<t'=t'(t)<t$,
satisfying 
$$x_{j_0}(t')-x_{j_0-1}(t')-{\textstyle{\sigma_0\over 2}} (t-t')=2y_0.$$
Indeed, for $t$ large enough, $x_{j_0}(t)-x_{j_0-1}(t)\ge {\sigma_0
\over 2} t \ge 2 y_0$,
and  $x_{j_0}(0)-x_{j_0-1}(0)-{\sigma_0\over 2}  t<0<2 y_0 $.
Then,
\bea\int \psi(.-(x_{j_0-1}(t)+y_0)) u^2(t)
&\le& \int \psi(.-(x_{j_0-1}(t')+{\textstyle{\sigma_0\over 2}}
(t-t')+y_0)) u^2(t')+C e^{-\gamma
y_0} \nonumber \\
&\le& \int \psi(.-(x_{j_0}(t')-y_0)) u^2(t') + C e^{-\gamma_0 y_0}.
\label{monottt}
\eea
Let $\delta>0$ be arbitrary. By $L^2_{loc}$ convergence of
$\e(t,.+x_{j_0}(t))$ and the induction assumption, we have, for
fixed $y_0$,
$$\int_{x>x_{j_0}(t)+2 y_0} \e^2(t) \goto 0 \quad \hbox{as $t\goto +\infty$.}
$$
Therefore, by Proposition \ref{CONVAROUND},
there exists $T=T(\delta)>0$, such that
$\forall t>T$, $\forall y_0>y_0(\delta)$,
\be\label{monotttt}\Big|
\int \psi(.-(x_{j_0}(t)-y_0)) u^2(t)
- \Big( \int Q^2 \Big)\sum_{k=j_0}^N (c_k^{+\infty})^{5-p\over 2(p-1)} 
\Big| \le \delta.\ee
Moreover, since $t'(t)\goto +\infty$ as $t\goto +\infty$, by
possibly taking a larger $T(\delta)$, we also have
\be\label{monottttt}\Big|
\int \psi(.-(x_{j_0}(t')-y_0)) u^2(t')
- \Big( \int Q^2 \Big)\sum_{k=j_0}^N (c_k^{+\infty})^{5-p\over 2(p-1)} 
\Big| \le \delta,
\ee
and so 
\be\label{monott}
\Big|\int 
\psi(.-(x_{j_0}(t)-y_0)) u^2(t) - \int \psi(.-(x_{j_0}(t')-y_0)) u^2(t')
\Big|\le 2 \delta.
\ee
Thus, by (\ref{monottt}), we have
\be\label{67bis}
 \int \psi(.-(x_{j_0-1}(t)+y_0)) u^2(t) 
\le \int \psi(.-(x_{j_0}(t)-y_0)) u^2(t) + 2 \delta +C e^{-\gamma_0 y_0}.
\ee
Since $\psi(x)\ge 1/2$ for $x>0$, by the decay properties of $Q$ and
(\ref{67bis}), we obtain
\bee && \int_{x_{j_0-1}(t)+y_0<y<x_{j_0}(t)-y_0}
\e^2(t) \\ && \qquad \le  2 \Big( \int \psi(.-(x_{j_0-1}(t)+y_0)) u^2(t) 
-\int \psi(.-(x_{j_0}(t)-y_0)) u^2(t)\Big) +C e^{-\gamma_0 y_0} 
\\&& \qquad \le 4 \delta +C' e^{-\gamma_0 y_0}.
\eee
Thus,
$\int_{x>x_{j_0-1}(t)} \e^2(t) \goto 0$ as $t\goto +\infty$.

\medskip

Finally, we prove $\int_{x>c_1^0 t/10} \e^2(t)\goto 0$ as
$t\goto +\infty$.
Indeed, let $0<t'=t'(t)<t$ such that
$x_1(t')-{c_1^0\over 20} (t +t')=y_0$,
Then, for $\sup_{t\ge 0} |\e(t)|_{H^1}$ small enough,
\bee 
\int \psi\Big(x- {c_1^0\over 10} t\Big)u^2(t) &\le &
\int \psi\Big(x-\Big({c_1^0\over 10}t' +{c_1^0\over 20} (t-t')\Big)\Big)
u^2(t')+C e^{-\gamma_0 y_0} \\
&\le & \int \psi(x-(x_1(t')-y_0)) u^2(t') +C e^{-\gamma_0 y_0}.
\eee
Arguing as before, this is enough to conclude the proof.

\bigskip

{\bf Proof of Corollary 1.} \quad 
Note first that 
\be\label{solprofile}
\Big | U^{(N)}(\,.\, ; c_j^0,-y_j)
-\sum_{j=1}^N Q_{c_j^0}(.-y_j) \Big |_{H^1}
\goto 0 \quad \hbox{as $\inf(y_{j+1}-y_j) \goto +\infty.$}
\ee
For $\gamma_0$, $A_0$, $L_0$ and $\alpha_0$ as in the statement of
Theorem 1, let $\alpha<\alpha_0$, $L>L_0$ be such that
$A_0\Big ( \alpha + e^{-\gamma_0 L} \Big ) < \delta_1/2$ and
\be\label{defofL}
\Big | U^{(N)}(\,.\, ; c_j^0,-y_j)
-\sum_{j=1}^N Q_{c_j^0}(.-y_j) \Big |_{H^1}
\le \delta_1/2, \quad \hbox{for $y_{j+1}-y_j>L$.}
\ee
Let $v(t,x)= U^{(N)}(x ; c_j^0,-(x_j^0+c_j^0 t))$ 
be an $N$-soliton solution. Let 
$T>0$ be such that
\be\label{by}
\forall t\ge T_1,\quad 
\Big |v(t)-\sum_{j=1}^N Q_{c_j^0}(.-(x_j^0+c_j^0 t))\Big |_{H^1}
\le \alpha/2,
\ee
and $\forall j,$ 
$x_{j+1}^0 + c_{j+1}^0 T \ge x_{j}^0 + c_{j}^0 T + 2L.$

By continuous dependence of the solution  of (\ref{kdvp})
with respect to the initial data (see \cite{KPV}), there
exists $\alpha_1>0$ such that if 
$|u(0)-v(0)|_{H^1}\le \alpha_1$, then
$|u(T)-v(T)|_{H^1}\le \alpha/2$.
Therefore, by (\ref{by}) 
$$\Big |u(T)-\sum_{j=1}^N Q_{c_j^0}(.-(x_j^0+c_j^0 T))\Big |_{H^1}
\le \alpha.$$
Thus, by Theorem 1 (i), there exists $x_j(t)$, for all $t\ge T$
such that
$$
\forall t\ge T,\quad
\Big |u(t)-\sum_{j=1}^N Q_{c_j^0}(.- x_j(t))\Big |_{H^1}
\le A_0 \Big ( \alpha + e^{-\gamma_0 L} \Big ) < \delta_1/2.
$$
Moreover, $x_{j+1}(t)>x_j(t)+L$.
Together with (\ref{defofL}), this gives the stability result.

Finally, Theorem 1 (ii) and (\ref{solprofile})
prove the asymptotic stability of the family of $N$-solitons.

\section*{Appendix : Modulation of a solution close to the sum of
$N$ solitons}

In this appendix, we prove the following lemma.

Let $0<c_1^0<\ldots<c_N^0$, 
$\sigma_0={1\over 2}\min(c_1^0,c_2^0-c_1^0,c_3^0-c_2^0,\ldots,c_N^0-c_{N-1}^0)
$. For $\alpha, L>0$, we consider the neighborhood 
of size $\alpha$ of the superposition of $N$ solitons
of speed $c_j^0$, located at a distance larger than $L$
\be\label{33}
{\cal U}(\alpha,L)=\big\{u\in H^1(\R);
\inf_{x_j>x_{j-1}+L }
\Big|u- \sum_{j=1}^N Q_{c_j^0}(.-x_j)
\Big|_{H^1}\le \alpha \Big\}.
\ee
(Note that functions in ${\cal U}(\alpha,L)$
have no time dependency.)

\begin{lemma}[Choice of the modulation parameters]\label{MODULATION}
 There exists $\alpha_1>0$, $L_1>0$ and unique
$C^1$ functions $(c_j,x_j):{\cal U}(\alpha_1,L_1)
\goto (0,+\infty)\times \R ,$ 
such that if $u\in {\cal U}(\alpha_1,L_1)$,
and 
\be\label{32}
\e(x)
=u(x)- \sum_{j=1}^N Q_{c_j}(.-x_j),
\ee
then
\be\label{orthoA}
\int Q_{c_i}(x-x_i) \e(x)dx=
\int (Q_{c_i})_x(x-x_i) \e(x)dx=0.
\ee
Moreover, there exists $K_1>0$ such that if $u\in
{\cal U}(\alpha,L)$, with $0<\alpha<\alpha_1$, $L>L_1$,
then
\be\label{smallnessA}
|\e|_{H^1}+\sum_{j=1}^N |c_j-c_j^0|
\le K_1 \alpha,
\quad x_j>x_{j-1} + L -K_1\alpha.
\ee
\end{lemma}

{\bf Proof.} Let $u \in {\cal U}(\alpha,L)$. It is clear
that for $\alpha$ small enough and $L$ large enough, the infimum
 $$\inf_{x_j\in \R}\Big|u-\sum_{j=1}^N Q_{c_j^0}(.-x_j)
\Big|_{H^1}$$
is attained for $(x_j)$ satisfying $x_j>x_{j-1}+L-C \alpha$,
for some constant $C>0$ independent of $L$ and $\alpha$.
By using standard arguments involving the implicit function theorem,
there exists $\alpha_1,L_1>0$ such that there exists unique $C^1$ 
functions $(r_j): 
{\cal U}(\alpha_1,L_1)\goto \R$, such that for all
$ u\in {\cal U}(\alpha,L)$, for $0<\alpha<\alpha_1$, $L>L_1$, we have
$$\Big|u-\sum_{j=1}^N Q_{c_j^0}(.-r_j(u))\Big|_{H^1}
=\inf_{x_j\in \R}\Big|u- \sum_{j=1}^N Q_{c_j^0}(.-x_j)
\Big|_{H^1}
\le \alpha.$$
Moreover, $r_j(u)-r_{j-1}(u)>L - C \alpha.$

\medskip

For some $c_j$,  $y_j$, $u\in H^1(\R)$, let 
$$Q_{c_j,y_j}(x)=Q_{c_j}(x-r_j(u)-y_j),\quad 
\e(x)=u(x)-\sum_{j=1}^N Q_{c_j,y_j}(x).$$
Define the following functionals:
$$\rho^{1,j}(c_1,\ldots,c_N,y_1,\ldots,y_N,u)=\int Q_{c_j,y_j}(x) \e(x) dx,
$$
$$\rho^{2,j}(c_1,\ldots,c_N,y_1,\ldots,y_N,u)=\int (Q_{c_j,y_j})_x(x) \e(x) dx,
$$
and $\rho=(\rho^{1,1},\rho^{2,1},\ldots,\rho^{1,N},\rho^{2,N}).$
Let $M_0=(c_1^0,\ldots,c_N^0,0,\ldots,0,\sum_{j=1}^N 
Q_{c_j^0,0}).$
We claim the following.

\medskip

{\sl Claim.} (i) $\forall j$,
$$
{\partial \rho^{1,j} \over \partial c_j}(M_0)=
-{5-p\over 4(p-1)}  (c_j^0)^{7-3p \over 2(p-1)} \int Q^2,\quad 
{\partial \rho^{1,j} \over \partial y_j} (M_0)=0,
$$
$$
{\partial \rho^{2,j} \over \partial c_j} (M_0)=0,
\quad  {\partial \rho^{2,j} \over \partial y_j} (M_0) 
= (c_j^0)^{p+3\over 2(p-1)} \int Q_x^2.
$$

(ii) $\forall j\not =k$,
$$\Big| {\partial \rho^{1,j} \over \partial c_k}(M_0)\Big|
+\Big| {\partial \rho^{1,j} \over \partial y_k}(M_0)\Big|
+\Big| {\partial \rho^{2,j} \over \partial c_k}(M_0)\Big|
+\Big| {\partial \rho^{2,j} \over \partial y_k}(M_0)\Big|\le
C e^{-\sqrt{\sigma_0} L/2}.$$

\medskip

{\sl Proof of the claim.}\quad 
Since 
$${\partial  Q_{c_j,y_j}\over \partial c_j}|_{(c^0_j,0)}
={1\over 2c_j^0} \left({2\over p-1} Q_{c_j^0,0}
 +(x-r_j) (Q_{c_j^0,0})_x\right),
\quad{\partial  Q_{c_j,y_j}\over \partial y_j}|_{(c^0_j,0)}
=-(Q_{c_j^0,0})_x,$$
we have by direct calculations:
\bee
{\partial \rho^{1,j} \over \partial c_j} (M_0) &=&  -\int Q_{c_j^0,0}
{\partial  Q_{c_j,y_j}\over \partial c_j}|_{(c^0_j,0)} 
 =- {1\over 2c_j^0} \int Q_{c_j^0,0}
    \left({2\over p-1} Q_{c_j^0,0} +(x-r_j) (Q_{c_j^0,0})_x\right) \\
&=& -{1\over 2} (c_j^0)^{7-3p \over 2(p-1)} \int Q  
\left({2\over p-1} Q +x Q_x\right) = -{1\over 2} (c_j^0)^{7-3p \over 2(p-1)}
{5-p\over 2(p-1)} \int Q^2,
\eee
by change of variable and integration by parts.
For  $j\not =k$
\bee
\left|{\partial \rho^{1,j} \over \partial c_k} (M_0)\right| 
&=&{1\over 2 c_k^0} \left| \int Q_{c_j^0,0}
\left({2\over p-1} Q_{c_k^0,0} +(x-r_k) (Q_{c_k^0,0})_x\right)\right|
\\ &\le& C \int e^{-\sqrt{\sigma_0}(|x-r_j|+|x-r_k|)} dx
\le  e^{-\sqrt{\sigma_0} |r_j-r_k|/2 }\le C e^{-\sqrt{\sigma_0} L/2}.
\eee
The rest is done in a similar way, using
$\int Q Q_x=0$, and $\int (Q_c)_x^2= c^{p+3\over 2(p-1)} \int Q_x^2$.

\medskip

It follows that $\nabla \rho(M_0)=D+P$,
where $D$ is a diagonal matrix with nonzero coefficients 
of order one on the 
diagonal, and $\| P\|\le C e^{-\sqrt{\sigma}_0 L/2}$.
Therefore, for $L$ large enough, the absolute value of
the Jacobian of $\rho$
at $M_0$ is larger 
than a positive constant depending only on the $c_j^0$.
Thus, by the implicit function theorem, by possibly taking a smaller
$\alpha_1$, there exists  $C^1$ functions $(c_j,y_j)$ of 
$u\in {\cal U}(\alpha_1,L_1)$ in  
a neighborhood of $(c_1^0,\ldots,c_N^0,0,\ldots,0)$ such that
$\rho(c_1,\ldots,c_N,y_1,\ldots,y_N,u)=0$.
Moreover, for some constant $K_1>0$, if
$u\in {\cal U}(\alpha,L_1)$, where $0<\alpha<\alpha_1$, then
$$\sum_{j=1}^N |c_j-c_j^0|+\sum_{j=1}^N |y_j| \le K_1 \alpha.$$
The fact that $|\e|_{H^1}\le K_1 \alpha$ then follows from
its definition. 
Finally, we choose $x_j(u)=r_j(u)+y_j(u)$.

\bigskip

Yvan Martel \\ Universit\'e de Cergy--Pontoise,
D\'epartement de Math\'ematiques \\ 2, av. Adolphe Chauvin,
95302 Cergy--Pontoise, France

\bigskip

Frank Merle \\
Institut Universitaire de France

and
\\ Universit\'e de Cergy--Pontoise,
 D\' epartement de Math\'ematiques \\
(same address)

\bigskip

Tai-Peng Tsai \\
Courant Institute,
New York University \\
251 Mercer Street, New York, NY 10012, USA


\end{document}